\documentclass[11pt]{amsart}
\usepackage{amsmath}
\usepackage{amsfonts}

\usepackage{amscd}
\usepackage{amsthm}
\usepackage{amssymb} \usepackage{latexsym}
\usepackage{eufrak}
\usepackage{euscript}
\usepackage{epsfig}
\usepackage{graphics}
\usepackage{array}
\usepackage{enumerate}
\usepackage{dsfont}
\usepackage{color}
\usepackage{wasysym}
\usepackage{hyperref}
\usepackage{pdfsync}

\newcommand{\bel}[1]{\begin{equation}\label{#1}}

\newcommand{\be}{\begin{equation}}

\newcommand{\ba}{\begin{eqnarray}}
\newcommand{\ea}{\end{eqnarray}}

\newcommand{\qe}{\end{equation}}

\newcommand{\Hmm}[1]{\leavevmode{\marginpar{\tiny%
$\hbox to 0mm{\hspace*{-0.5mm}$\leftarrow$\hss}%
\vcenter{\vrule depth 0.1mm height 0.1mm width \the\marginparwidth}%
\hbox to
0mm{\hss$\rightarrow$\hspace*{-0.5mm}}$\\\relax\raggedright #1}}}

\theoremstyle{theorem}
\newtheorem{thm}{Theorem}[section]

\theoremstyle{example}

\theoremstyle{corollary}

\theoremstyle{lemma}
\newtheorem{lem}[thm]{Lemma}
\theoremstyle{definition}
\newtheorem{defi}[thm]{Definition}
\theoremstyle{proof}
\theoremstyle{remark}
\newtheorem{rem}[thm]{Remark}

\begin{document}

\title{Curvature Estimate on the finite graph with large girth}
\author{Yijin Gao}
\email{yjgao@ruc.eud.cn}
\address{Department of Mathematics,Information School,
Renmin University of China,
Beijing 100872, China
}

\begin{abstract}The CD inequalities and CDE inequalities are useful in the  estimate of curvature on graphs. This article is based on the ufinite graph with large girth, and finally concludes some curvature estimate in CD and CDE.
\end{abstract}
\maketitle

\section{Introduction}
Ricci curvature on graphs is an important aspect in the study of graphs and geometry.Some related work has been done in {\cite{6}}{\cite{7}}.
Recently,there is some work has been done about graphs with large girth in {\cite{2}}.So we want to do some research about the curvature estimate under the special conditions-when the girth of a finite graph is large than 5.Maybe the conclusion is beautiful .

Moreover,the graph referred is unweighted,the corresponding Laplacian is called unweighted normalized(i.e.$\mu=1$ on E and $m=\mu$ on V){\cite {1}}{\cite{3}}{\cite{4}}.
Let $(V,E)$ be a undirected graph with the set of vertices $V$ and the set of edges $E$,i.e. two-elements subsets in $V$.The graph is called simple if there is no self-loops and multiple edges.The graph is called locally finite,if the combinatorial degree $d_x<\infty$ for any $x \in V$.We say a vertex $x$ is a pending vertex if $d_x=1$.For any subsets $A,B\subset V$,we denote by $E(A,B):={{x,y}\in E: x\in A ,y\in B}$ the set of edges between $A$ and $B$.For vertices $x$ and $x$ and $y$,a walk from $x$ to $y$ is a sequence of vertices ${x_i}$such that $$x=x_o\sim{x_1}\sim...\sim{x_k}=y,$$
where $k$ is called the length of the walk.A graph is said to be connected if for any $x,y \in V$ there is a walk from $x$ to $y$.In this paper,we only consider undirected,connected,locally finite simple graphs.

\indent This paper gives an estimate of the curvature  on finite graph with  girth large than 5. \\
\indent The paper is organized into four parts:\\
\indent Chapter $1$ is the introduction of the finite graph, the girth,the Laplacians , CD inequalities and  CDE inequalities on it.\\
\indent Chapter $2$ introduces some easy conclusions about the calculation of the operator which is helpful for chapter $3$.\\
\indent Chapter $3$ is the main conclusion of this thesis which includes curvature estimate about CD inequality and CDE inequality.

Acknowledgement:The author is grateful for the help and suggestions from his advisor Pro.Yong Lin.

\section{GRAPHS,GIRTH, LAPLACIANS CD INEQUALITIES AND CDE INEQUALITIES}

\indent Given a graph $G=(V,E)$, for an $x\in V$, if there exists another $y\in V$ that satisfies $(x,y)\in E$, we call them are neighbors, and written as $x\sim y$. If there exists an $x\in V$ satisfying $(x,x)\in E$, we call it a self-loop. \\
\indent Now we will introduce some basic definitions and theorems before we get the main results.
\begin{defi} {\rm (locally finite graph)}\ \ We call a graph G is a locally finite graph if for any $x\in V$, it satisfies $\#\{y\in V|y\sim x\}<\infty$. Moreover, it is called connected if there exists a sequence $\{x_i\}_{i=0}^n$ satisfying: $x=x_0\sim x_1\sim \cdots \sim x_n=y$.\end{defi}
\begin{defi} {\rm (Laplacians on locally finite graphs)}\ \ On a locally finite graph $G=(V,E,\mu,m)$ the Laplacian has a form as follows: $$\triangle f(x)=\frac{1}{m(x)}\sum_{y\in V}\mu_{xy}(f(y)-f(x)),\quad \forall f\in C_0(V).$$\end{defi}
\begin{defi} {\rm (gradient operator $\Gamma$)}\ \ The operator $\Gamma$ is defined as follows:$$\Gamma (f,g)(x)=\frac{1}{2}(\triangle(fg)-f\triangle g-g\triangle f)(x)$$.\end{defi}
\indent Always we write $\Gamma (f,f)$ as $\Gamma(f)$.
\begin{defi} {\rm (gradient operator $\Gamma_{i}$)}\ \ The operator $\Gamma_{i}$ is defined as follows:$$\Gamma_{0}(f,g)=fg$$\\ $$\Gamma_{i+1}(f,g)=\frac{1}{2}(\triangle(\Gamma_{i}(f,g))-\Gamma_{i}(f,\triangle g)-\Gamma_{i}(\triangle f,g))$$.\end{defi}
\indent Also we have $\Gamma_2(f)=\Gamma_2(f,f)=\frac{1}{2}\triangle\Gamma(f)-\Gamma(f,\triangle f)$.
\begin{defi} {\rm (the girth of a graph)}\ \ The girth of a verth $x$ in $(V,E)$,denoted by $Gir(x)$,is defined to the minimal length of cycles passing through $x$.(If there is no cycle passing through $x$,define $Gir(x)=\infty$.)The girth of a graph is defined as $\inf_{x \in V} Gir(x)$. \end{defi}
\begin{defi}{\rm ($CD(K,n)$ condition)} \ \ We call a graph satisfies $CD(K,n)$ condition if for any $x \in V$,we have$$\Gamma_2(f)(x)\ge \frac{1}{n}(\triangle f)^2(x)+K \Gamma(f)(x).\quad K\in \mathbb{R}.$$\end{defi}
\begin{defi} {\rm ($CDE(K,n)$ condition)}\ \ Let $f:V\rightarrow\mathbb{R}^+$ satisfy $f(x)>0$, $\triangle f(x)<0$. We call a graph satisfies $CDE(x,K,n)$ condition if for any $x\in V$, we have$$\Gamma_2(f)(x)-\Gamma\left(f,\frac{\Gamma(f)}{f}\right)(x)\geqslant\frac{1}{n}(\triangle f)(x)^2+K\Gamma(f)(x).\quad K\in \mathbb{R}.$$\end{defi}

\section{MAIN RESULTS}

\begin{rem} In this section,we use the definition in part 2.They are the work from {\cite {7}}.\\
\end{rem}

\begin{lem} $$\Gamma(f)(x)=\frac{1}{2}\mu_x\sum_{y \in N_x} f(x,y)^2.$$ \\
Here we define $N_x={y \in V:xy \in E}$ and $d_x=|N_x|.$For notational simplicity we work with $\mu_x=\frac{1}{d_x}.$
\end{lem}
{\bf Proof}
\indent We have
\begin{equation*}
\left.
\begin{aligned}
&\Gamma (f)(x)=\frac{1}{2}\triangle(f^2)(x)-f(x)(\triangle f)(x)\\
&=\frac{1}{2}\mu_x \sum_{y \in N_x}(f^2)(x,y)-f(x)\mu_x \sum_{y \in N_x} f(x,y)\\
&=\frac{1}{2}\mu_x\sum_{y \in N_x}(f(x,y)(f(y)+f(x))-2f(x,y)f(x))\\
&=\frac{1}{2}\mu_x\sum_{y \in N_x}f(x,y)^2.\\
\end{aligned}
\right.
\end{equation*}

\begin{lem} $$\Gamma_2(f)(x)=\frac{1}{2}((\triangle f)^2(x)+\mu_x \sum_{y \in N_x}\mu_y\sum_{z \in N_y} (f(y,z)^2-\frac{1}{2}f(x,z)^2)).$$
{\bf Proof}
\indent We have
\begin{equation*}
\left.
\begin{aligned}
&\triangle (\Gamma(f))(x)=\mu_x\sum_{y \in N_x} \Gamma(f)(x,y)
&=\mu_x \sum_{y \in N_x}\frac{1}{2}\mu_y\sum_{z \in N_y}(f(y,z)^2-f(x,y)^2)\\
\end{aligned}
\right.
\end{equation*}
\indent and
\begin{equation*}
\left.
\begin{aligned}
&\Gamma(f,\triangle f)(x)=\frac{1}{2}(\triangle(f\cdot \triangle f)(x)-f(x)\cdot(\triangle^2 f)(x)-(\triangle f)^2(x))\\
&=-\frac{1}{2}(\triangle f)^2(x)+\frac{1}{2}\mu_x \sum_{y \in N_x}((f \triangle f)(x,y)-f(x)(\triangle f)(x,y))\\
&=-\frac{1}{2}(\triangle f)^2(x)+\frac{1}{2}\mu_x \sum_{y \in N_x}f(x,y)(\triangle f)(y)\\
&=-\frac{1}{2}(\triangle f)^2(x)+\frac{1}{2}\mu_x\sum_{y \in N_x}f(x,y)\mu_y \sum_{z \in N_y}f(y,z)\\
\end{aligned}
\right.
\end{equation*}

\indent thus
\begin{equation*}
\left.
\begin{aligned}
&\Gamma_2(f)(x)=\frac{1}{2}\triangle (\Gamma(f))(x)-\Gamma(f,\triangle f)(x)\\
&=\frac{1}{2}(\triangle f)^2(x)+\frac{1}{2}\mu_x\sum_{y \in N_x}\mu_y \sum_{z \in N_y}(\frac{1}{2}f(y,z)^2-\frac{1}{2}f(x,y)^2-f(x,y)f(y,z))\\
&=\frac{1}{2}(\triangle f)^2+\frac{1}{2}\mu_y\sum_{y \in N_x} \mu_y \sum_{z \in N_y}(f(y,z)^2-\frac{1}{2}f(x,z)^2)\\
\end{aligned}
\right.
\end{equation*}

\end{lem}

\section{BASIC CONCLUSION}
\begin{rem} In this section,we will give the curvature estimate about finite graphs with girth larger than 5 in the CD inequality .\\
\end{rem}

\begin{thm}
\indent we have a finite with girth larger than 5 and we concern the fixed point $x$,assume the neighbood of $x$ are $y_1,y_2,...,y_n$,so just from the definition above ,we have $d_x=n$.Also,the neighbood of $y_1$ is $z_{11},z_{12},...,z_{1k_1-1}$,the neighbood of $y_2$ is $z_{21},z_{22},...,z_{2k_2-1}$...,the neighbood of $y_n$ is $z_{n1},z_{n2},...,z_{nk_n-1}$.So $d_{y_{i}}=k_{i},i=1,2,...,n$.Let $k=min\{\frac{2-k_1}{k_1},\frac{2-k_2}{k_n},...\frac{2-k_n}{k_n}\}$. Then we have the conclusion that the graph satisfies $CD(k,2,x)$.\\
{\bf Proof}
\indent From the theorem above,we know that
$$\Gamma_2(f)(x)=\frac{1}{2}(\triangle f)^2(x)+\frac{1}{2}\mu_y\sum_{y \in N_x}\mu_y\sum_{z \in N_y}(f(y,z)^2-\frac{1}{2}f(x,z)^2).$$
\indent So with the definition of $CD(k,2,x)$ inequality ,we need to prove the following inequality
$$\frac{1}{2}(\triangle f)^2+\frac{1}{2}\mu_x\sum_{y \in N_x} \mu_y \sum_{z \in N_y}(f(y,z)^2-\frac{1}{2}f(x,z)^2))\ge \frac{1}{2}(\triangle f)^2+k\Gamma(f).$$\\
\indent Which is the same to:
$$\frac{1}{2}\mu_x\sum_{y \in N_x}\mu_y\sum_{z \in N_y}(f(y,z)^2-\frac{1}{2}f(x,z)^2))\ge k\cdot \frac{1}{2}\mu_x\sum_{y \in N_x} f(x,y)^2.$$\\
\indent Here $\mu_x =\frac{1}{n},\mu_{y_{i}}=\frac{1}{k_i} (i=1,2,..,n)$,and without loss of generality we can aasume $f(x)=0$.\\
\indent So the inequality is equal to the following:
$$\sum_{y \in N_x} \mu_y \sum_{z \in N_y}(f(z)-f(y))^2-\frac{1}{2}f(z)^2) \ge k\cdot \sum_{y \in N_x}f(y)^2.$$ \\
\indent Firstly,we concern the $y_1$ part in the left inequality.\\
The sum is :
\begin{equation*}
\left.
\begin{aligned}
&\frac{1}{k_1}((f(z_{11}-f(y_1))^2-\frac{1}{2}f(z_{11})^2+(f(z_{12})-f(y_1))^2+(f(z_{21})-f(y_1))^2-\frac{1}{2}f(z_{12})^2+...)\\
&=\frac{1}{k_1}(f(z_{11})^2+f(y_1)^2-2f(z_{11})f(y_1)-\frac{1}{2}f(z_{11})^2+...f(y_1)^2)\\
\end{aligned}
\right.
\end{equation*}
\indent We use the knowlegde of quadratic function,take $f(z_{1i})=2f(y_1)$,for $i=1,2,...,k_1-1$,so we can get the minimal sum of the function.\\
The sum about $y_1$ is equal to:
$$\frac{1}{k_1}(-f(y_1)^2\cdot (k_1-1)+f(y_1)^2)=\frac{1}{k_1}f(y_1)^2(2-k_1)=-\frac{k_1-2}{k_1}f(y_1)^2.$$
\indent And the situation is all the same to $y_2,y_3,...,y_n$.At last we have the following inequality:
$$\frac{2-k_1}{k_1}f(y_1)^2+\frac{2-k_2}{k_2}f(y_2)^2+...+\frac{2-k_n}{k_n}f(y_n)^2\ge kf(y_1)^2+kf(y_2)^2+...kf(y_n)^2.$$
\indent From the assumation of the theorem,we know that $\frac{2-k_i}{k_i} \ge k $ for $i=1,2,...,n.$ Also $f(y_i)^2 \ge 0.$So the above inequality is correct.\\we end the proof.
\end{thm}

\begin{thm}
\indent Suppose $G$ is a graph with girth large than 5 and let $f:V\to R^{+}$ satisfy $f(x) >0,\triangle f(x) <0$.Then we have the conclusion that the graph satisfy $CDE(2,-\frac{n}{2}-1)$\\
{\bf Proof}
\indent Because the $CDE$ condition satisfies that:
$$\Gamma_2(f)-\Gamma(f,\frac{\Gamma(f)}{f}) \ge \frac{1}{m}(\triangle f)^2+k\cdot \Gamma(f).$$
\indent From the above proposition ,we hava that:
$$\Gamma_2(f)(x)=\frac{1}{2}((\triangle f)^2(x)+\mu_x\sum_{y \in N_x}\mu_y \sum_{z \in N_y}((f(y,z)^2-\frac{1}{2}f(x,z)^2)).$$
\indent So the most important work for us is to simplify the $\Gamma(f,\frac{\Gamma(f)}{f})$.
\begin{equation*}
\left.
\begin{aligned}
&\Gamma(f,\frac{\Gamma(f)}{f})(x)=\frac{1}{2}(\triangle(\Gamma(f))(x)-\triangle(\frac{\Gamma(f)}{f})-\triangle f \cdot \Gamma(f))\\
&=\frac{1}{2}\triangle (\Gamma(f))(x)-\frac{1}{2}\triangle (\frac{\Gamma(f)}{f})-\frac{1}{2}\triangle f \cdot \Gamma(f)\\
&=I_1-I_2-I_3\\
\end{aligned}
\right.
\end{equation*}

\indent Firstly:
\begin{equation*}
\left.
\begin{aligned}
&I_1=\frac{1}{2}\triangle(\Gamma(f))(x)=\frac{1}{2}\mu_x \sum_{y \in N_x}\Gamma(f)(x,y)\\
&=\frac{1}{2}\mu_x \sum_{y \in N_x}\frac{1}{2} \mu_y \sum_{z \in N_y} (f(y,z)^2-f(x,y)^2)\\
\end{aligned}
\right.
\end{equation*}

\indent Secondly:
\begin{equation*}
\left.
\begin{aligned}
&I_2=\frac{1}{2}\triangle  (\frac{\Gamma(f)}{f})=\frac{1}{2}\frac{1}{n}\sum_{y \in N_x}\frac{\Gamma(f)}{f}(x,y)\\
&=\frac{1}{2n}\sum_{y \in N_x}\frac{\Gamma(f)(y)}{f(y)}-\frac{1}{2n}\sum_{y \in N_x}\frac{\Gamma(f)(x)}{f(x)}\\
&=\frac{1}{2n}\sum_{y \in N_x}\frac{1}{2f(y)}\mu_y \sum_{z \in N_y} f(y,z)^2-\frac{1}{2n}\frac{1}{2}\sum_{y \in N_x} f(x,y)^2\\
&=\frac{1}{4n}\sum_{y \in N_x} \mu_y \sum_{z \in N_y} \frac{f(y,z)^2}{f(y)}-\frac{1}{4n}\sum_{y \in N_x} f(x,y)^2\\
\end{aligned}
\right.
\end{equation*}

\indent Thirdly:
$$I_3=\frac{1}{2}\triangle f(x)\Gamma(f)(x)=\frac{1}{2n}\sum_{y \in N_x}f(x,y) \frac{1}{2n}\sum_{y \in N_x} f(x,y)^2.$$

\indent To simplify the question,we assume that $f(x)=1;f(x,y_i)=v_i,i=1,2,..n.$\\
so combine these together,we have:

\noindent$\displaystyle
k\ge \frac{1}{\frac{1}{2n}\sum_{y \in N_x}{v_i}^2}(\frac{1}{2}(\triangle f)^2+
\frac{1}{2n}\sum_{y \in N_x}\mu_y \sum_{z \in N_y}
(f(y,z)^2-\frac{1}{2}f(x,z)^2)+\\
\frac{1}{4n}\sum_{y \in N_x}\mu_y \sum_{z \in N_y} \frac{f(y,z)^2}{f(y)}-\frac{1}{4n}
\sum_{y \in N_x}f(x,y)^2+\frac{1}{2n}\sum_{y \in N_x} f(x,y)\frac{1}{2n}f(x,y)^2-\frac{1}{2n^2}(\sum_{y \in N_x}f(x,y))^2)\\
=\frac{1}{2n\sum_{y \in N_x}{v_i}^2}(\frac{1}{2n}\sum_{y \in N_x} \mu_y \sum_{z \in N_y}((f(y,z)^2-\frac{1}{2}f(x,z)^2-\frac{1}{2}f(y,z)^2\\
+\frac{1}{2}f(x,y)^2+\frac{1}{2}\frac{f(y,z)^2}{f(y)}+\frac{1}{2n}\sum_{y \in N_x} {v_i}^2 \frac{1}{2n}\sum_{y \in N_x}v_i\\
-\frac{1}{4n}\sum_{y \in N_x} {v_i}^2)$

\indent Absoulte we have :
\begin{equation*}
\left.
\begin{aligned}
&f(y,z)^2-\frac{1}{2}f(x,z)^2-\frac{1}{2}f(y,z)^2+\frac{1}{2}f(x,y)^2+\frac{f(y,z)^2}{2f(y)}\\
&=\frac{f(z)^2}{f(y)}-f(y)f(z)+f(y)^2-\frac{1}{2}f(y)
\end{aligned}
\right.
\end{equation*}
\indent They are minimized when $f(z)=f(y)^2$,whence the sum is $$-\frac{1}{2}f(y)^3+f(y)^2-\frac{1}{2}f(y)=-\frac{1}{2}\sum{v_i}^3-\frac{1}{2}\sum{v_i}^2.$$
\indent So
$$k \ge -1+\frac{1}{\sum{v_i}^2}(\frac{1}{2n}\sum{v_i}^2\sum v_i -\frac{1}{2}\sum{v_i}^3.$$
\indent Because
$$f(y_i)>0,$$
\indent so
$$v_i=f(y_i)-f(x)=f(y_i)-1>-1.$$
\indent Also according to $\triangle f(x)=\frac{1}{n}\sum_{y \in N_x}f(x,y_i) <0$,we have $\sum v_i <0.$\\
we can assume $v_1<v_2<...<v_n$,so $v_n<n-1.$

\indent At last we need only to prove:
$$\frac{\sum v_i \sum {v_i}^2-n\sum {v_i}^3}{2n\sum {v_i}^2}>-\frac{n}{2}.$$
$$\sum v_i >-n \longrightarrow \frac{\sum v_i \sum{v_i}^2}{2n\sum {v_i}^2}>-\frac{n}{2n}=-\frac{1}{2}.$$
$$\sum {v_i}^3=\sum{v_i}^2v_i <\sum {v_i}^2 v_n <\sum {v_i}^2 (n-1) \longrightarrow \frac{-n\sum {v_i}^3}{2n\sum{v_i}^2}>-\frac{n(n-1)}{2n}=-\frac{n}{2}+\frac{1}{2}.$$

\indent we end the proof.

\end{thm}

\end{document}